\documentclass[a4paper,10pt,leqno,twoside]{tobart}
\usepackage[utf8]{inputenc}
\usepackage[english]{babel}
\usepackage{inputenc, amsmath, amssymb, latexsym,
  epic, epsfig, rotating, fancyhdr, amsthm, pifont}
\usepackage{paralist}
\usepackage{amsmath}
\usepackage{verbatim} 

\newcommand{\ld}{\ensuremath{,\ldots,}}
\newcommand{\ssq}{\ensuremath{\subseteq}}
\newcommand{\smin}{\ensuremath{\setminus}}
\newcommand{\eps}{\ensuremath{\varepsilon}}

\newcommand{\T}{\ensuremath{\mathbb{T}}}

\newcommand{\N}{\ensuremath{\mathbb{N}}} 
\newcommand{\R}{\ensuremath{\mathbb{R}}}
\newcommand{\Z}{\ensuremath{\mathbb{Z}}}
\newcommand{\Q}{\ensuremath{\mathbb{Q}}}

\newcommand{\Leb}{\ensuremath{\mathrm{Leb}}}

\newcommand{\diam}{\ensuremath{\mathrm{diam}}}

\newcommand{\kreis}{\ensuremath{\mathbb{T}^{1}}}

\newcommand{\alphlist}{\begin{list}{(\alph{enumi})}{\usecounter{enumi}\setlength{\parsep}{2pt}
      \setlength{\itemsep}{1pt} \setlength{\topsep}{5pt}
      \setlength{\partopsep}{3pt}}}

\newcommand{\arablist}{\begin{list}{(\arabic{enumi})}{\usecounter{enumi}\setlength{\parsep}{2pt}
          \setlength{\itemsep}{1pt} \setlength{\topsep}{5pt}
          \setlength{\partopsep}{3pt}}}

\newcommand{\romanlist}{\begin{list}{(\roman{enumi})}{\usecounter{enumi}\setlength{\parsep}{2pt}
              \setlength{\itemsep}{1pt} \setlength{\topsep}{5pt}
              \setlength{\partopsep}{3pt}}}

 \newcommand{\listend}{\end{list}}

\newcommand{\bulletlist}{\begin{list}{$\bullet$}{\setlength{\parsep}{2pt}
                \setlength{\itemsep}{1pt} \setlength{\topsep}{5pt}
                \setlength{\partopsep}{3pt}\setlength{\leftmargin}{15pt}}}

\newcommand{\roundqed}{{\hfill \Large $\circ$}}

\newcommand{\foot}{\footnote}

\newcommand{\ifolge}[1]{\ensuremath{(#1)_{i\in\mathbb{N}}}}

\newcommand{\nLim}{\ensuremath{\lim_{n\rightarrow\infty}}}

\newcommand{\ntel}{\ensuremath{\frac{1}{n}}}

\newcommand{\halb}{\ensuremath{\frac{1}{2}}}

\newtheoremstyle{tobthm}{3pt}{3pt}{\itshape}{0pt}{\bfseries}{.}{0.5eM}{}
\theoremstyle{tobthm}

\def\clap#1{\hbox to 0pt {\hss#1\hss}}

\newtheorem{definition}{Definition}[section]
\newtheorem{thm}[definition]{Theorem}

\newtheorem{lem}[definition]{Lemma}
\newtheorem{cor}[definition]{Corollary}  
\newtheorem{prop}[definition]{Proposition}

\newtheorem{claim}[definition]{Claim}

\newtheoremstyle{tobrem}{3pt}{3pt}{\normalfont}{0pt}{\bfseries}{.}{0.5em}{}
\theoremstyle{tobrem}

\newtheorem{rem}[definition]{Remark}    
\newtheorem{example}[definition]{Example}

\numberwithin{equation}{section}
\numberwithin{figure}{section}

\title{\Large\textsc{Dimensions of attractors in pinched skew products}}
\author{M.~Gr\"oger\thanks{Department of Mathematics, Universit\"at Bremen,
Germany. Email: {\tt groeger@math.uni-bremen.de}} \ and 
T.~J\"ager\thanks{Department of Mathematics, TU Dresden,
Germany. Email: {\tt Tobias.Oertel-Jaeger@tu-dresden.de}}}

\date{}

\newcommand{\set}[1]{\mathbb #1}
\newcommand{\abs}[1]{\left|#1\right|}
\newcommand{\Abs}[1]{\left\|#1\right\|}
\newcommand{\Acal}{\ensuremath{\mathcal{A}}}

\begin{document}

\setlength{\abovedisplayskip}{1.0ex}
\setlength{\abovedisplayshortskip}{0.8ex}

\setlength{\belowdisplayskip}{1.0ex}
\setlength{\belowdisplayshortskip}{0.8ex}

\maketitle

\abstract{We study dimensions of strange non-chaotic attractors and their
  associated physical measures in so-called pinched skew products, introduced by
  Grebogi and his coworkers in 1984.  Our main results are that the Hausdorff
  dimension, the pointwise dimension and the information dimension are all equal
  to one, although the box-counting dimension is known to be two.  The assertion
  concerning the pointwise dimension is deduced from the stronger result that
  the physical measure is rectifiable.  Our findings confirm a conjecture by
  Ding, Grebogi and Ott from 1989.}

\noindent
\section{Introduction}

\begin{figure}[h] 
\begin{center}
\epsfig{file=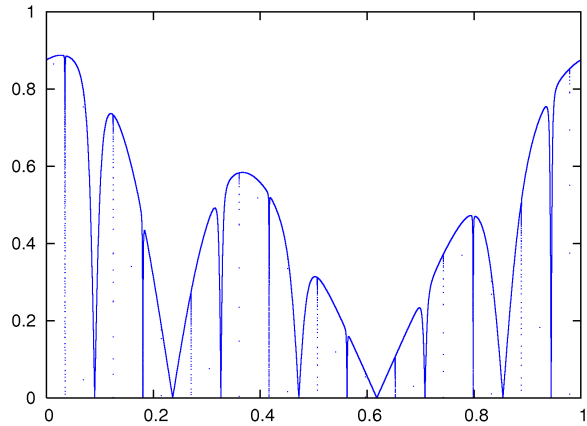, width=0.6\linewidth}
\end{center} {\caption{ \small Strange non-chaotic attractor in
    (\ref{e.upper-boundinggraph}) with $\kappa=3$ and $\rho$ the
    golden mean. \label{f.1}}}
\end{figure}

In \cite{grebogi/ott/pelikan/yorke:1984}, Grebogi and coworkers
introduced (a slight variation of) the system
\begin{equation} \label{e.gopy-map}
	F_\kappa \ : \ \kreis \times [0,1] \to \kreis \times [0,1] \quad , \quad F_\kappa(\theta,x)
	\ = \ (\theta+\rho\bmod 1,\tanh(\kappa x)\cdot\sin(\pi\theta))\ ,
\end{equation}
with $\rho\in\R\smin\Q$ and real parameter $\kappa>0$, as a simple model for the
existence of {\em strange non-chaotic attractors (SNA)}.\foot{The model studied
  by Grebogi {\em et al} was a four-to-one extension of (\ref{e.gopy-map}) with
  slightly different parametrisation.} Later, the term {\em `pinched skew
  products'} was coined by Glendinning \cite{glendinning:2002} for a general
class of systems sharing some essential properties of (\ref{e.gopy-map}). The
object which is called an SNA in the above system is the {\em upper bounding
  graph} $\varphi^+$ of the global attractor $\Acal:=\bigcap_{n\in\N}
F_\kappa^n(\kreis\times[0,1])$, which is given by
\begin{equation} \label{e.upper-boundinggraph}
  \varphi^+(\theta) \ := \ \sup\{x\in[0,1] \mid (\theta,x)\in\Acal\} \ .
\end{equation}
An illustration of this attractor is shown in Figure~\ref{f.1}.

Due to the monotonicity of the fibre maps $F_{\kappa,\theta}:x\mapsto
\tanh(\kappa x)\cdot\sin(\pi\theta)$, one can verify that the function
$\varphi^+$ satisfies
\begin{equation}\label{e.inv-graph}
F_{\kappa,\theta}(\varphi^+(\theta)) \ = \ \varphi^+(\theta+\rho\bmod 1) \ .
\end{equation}

Consequently, the corresponding point set
$\Phi^+:=\{(\theta,\varphi^+(\theta))\mid \theta\in\kreis\}$ is
$F_\kappa$-invariant.  Slightly abusing terminology, we will call both
$\varphi^+$ and $\Phi^+$ an {\em invariant graph}.  Keller showed in
\cite{Keller1996} that for $\kappa>2$ in (\ref{e.gopy-map}) the graph
$\varphi^+$ is $\Leb_{\kreis}$-almost surely strictly positive, its {\em
  Lyapunov exponent}
\begin{equation*}
  \lambda(\varphi^+) \ = \ \int\log F_{\kappa,\theta}'(\varphi^+(\theta)) \ d\theta
\end{equation*}
is strictly negative and $\varphi^+$ attracts $\Leb_{\kreis\times
  [0,1]}$-a.e.\ initial condition. Note that Birkhoff's Ergodic
Theorem implies that $\nLim \ntel\log
\left(F_{\kappa,\theta}^n\right)'(\varphi^+(\theta)) =
\lambda(\varphi^+)$ for $\Leb_{\kreis}$-a.e.\ $\theta\in\kreis$, where
we let $F^n_{\kappa,\theta} = F_{\kappa,\theta+(n-1)\rho\bmod 1} \circ \ldots
\circ F_{\kappa,\theta}$.

The findings in \cite{grebogi/ott/pelikan/yorke:1984} attracted substantial
interest in the theoretical physics community, and subsequently a large number
of numerical studies confirmed the widespread existence of SNA in
quasiperiodically forced systems and explored their behaviour and properties
(see \cite{prasad/negi/ramaswamy:2001,haro/puig:2006,jaeger:2006a} for an
overview and further references).  For a long time, however, rigorous results
remained rare, and even basic questions are still open nowadays. In particular,
this concerns the dimensions and fractal properties of SNA, which are still
mostly unknown even for the original example by Grebogi {\em et al}. A numerical
investigation was carried out in \cite{DingGrebogiOtt1989a}, and the results
indicated that the box dimension of the attractor is two, whereas the
information dimension should be one. For sufficiently large $\kappa$, the
conjecture on the box dimension was verified indirectly in \cite{jaeger:2007},
by showing that the topological closure of $\Phi^+$ is equal to the global
attractor $\Acal = \{(\theta,x)\mid 0\leq x \leq \varphi^+(\theta)\}$ and
therefore has positive two dimensional Lebesgue measure.

Our aim is to determine further dimensions of $\varphi^+$ and the associated
invariant measure $\mu_{\varphi^+}$, which is obtained by projecting the
Lebesgue measure on the base $\kreis$ onto $\Phi^+$. For the Hausdorff dimension
$D_H$ (see Section~\ref{Dimensions} for the definition), we have
\begin{thm}
\label{t.hausdorff-dim}
Suppose $\rho$ in (\ref{e.gopy-map}) is Diophantine and $\kappa$ is sufficiently
large. Then $D_H(\Phi^+) = 1$.
Furthermore, the one-dimensional Hausdorff measure of $\Phi^+$ is infinite.
\end{thm}
This statement is a special case of Corollary~\ref{c.final}, see
Section~\ref{Hausdorff}.  Here and in the results below, the largness condition
of $\kappa$ depends on the constants of the Diophantine condition on $\rho$.
\begin{rem}
  Our results in Section~\ref{Hausdorff} also allow to treat examples with a
  higher dimensional driving space, as given in Example
  \ref{example_higher_dimensional_case}. In these cases, the rotation on $\T^1$
  is replaced by a rotation on $\T^D$, and we obtain that the Hausdorff dimension
  of $\Phi^+$ is $D$. However, at least for sufficiently large $D$ the
  $D$-dimensional Hausdorff measure is finite, in contrast to the case $D=1$
  (Proposition \ref{proposition_Hausdorff_measure_finite}).  We believe that for
  these examples the $D$-dimensional Hausdorff measure is infinite only for
  $D=1$ and finite for all $D\geq 2$.
\end{rem}

In order to obtain information on the invariant measure $\mu_{\varphi^+}$, we
determine its pointwise dimension given by
\begin{equation*}
d_{\mu_{\varphi^+}}(\theta,x) \ = \ \lim_{\eps\to 0} \frac{\log
\mu_{\varphi^+}(B_\eps(\theta,x))}{\log \eps}.
\end{equation*}
A priori, it is not clear whether this limit exists, such that in general one
defines the upper and lower pointwise dimension by taking the limit superior and
inferior, respectively (see Section~\ref{Dimensions}). Furthermore, even if the
limit exists, it may depend on $(\theta,x)$. If the pointwise dimension exists
and is constant almost surely, the invariant measure is called {\em
exact dimensional}. It turns out that this is the case in the
situation considered here. 
In fact, we obtain the stronger result that $\mu_{\varphi^+}$ is a
rectifiabile measure, see Section \ref{rectifiability} and Theorem
\ref{theorem_rectifiability}, and this directly implies
\begin{thm}
\label{t.pointwise} Suppose $\rho$ in (\ref{e.gopy-map}) is Diophantine 
and $\kappa$ is sufficiently large. Then for $\mu_{\varphi^+}$-almost every
$(\theta,x)\in\kreis\times[0,1]$, we have $d_{\mu_{\varphi^+}}(\theta,x)=1$. In
particular, $\mu_{\varphi^+}$ is exact dimensional.
\end{thm}
For an exact dimensional measure $\mu$, it is known that the information
dimension $D_1$ (see again Section~\ref{Dimensions} for the definition)
coincides with the pointwise dimension. Hence, we obtain
\begin{cor}\label{t.info-dim}
Suppose $\rho$ in (\ref{e.gopy-map}) is Diophantine and $\kappa$ is sufficiently
large. Then $D_1(\mu_{\varphi^+})=1$.
\end{cor}
This confirms the conjecture made in \cite{DingGrebogiOtt1989a}.
Since the geometric mechanism for the creation of SNA in pinched skew
products is quite universal and can be found in similar form in other
types of systems, we expect our results to hold in further situations.
For example, this should be true for the SNA found in the Harper map,
which describes the projective action of quasiperiodic Schr\"odinger
cocycles, and for SNA in the quasiperiodically forced versions of the
logistic map and the Arnold circle map. On a technical level, these
systems are much more difficult to deal with, and for this reason we
refrain from extending our analysis beyond pinched skew products here.
Yet, combining our approach with the methods developed in
\cite{bjerkloev:2005a,bjerkloev:2005,jaeger:2009a} should allow to
produce similar results for the mentioned examples.  Apart from this,
progress has also been made recently concerning the existence of SNA
in quasiperiodically forced unimodal maps
\cite{bjerkloev2009sna,bjerkloev2012quasi,alseda2008attractors} .
Here, similar results may be expected as well, but it is much less
clear to what extend the presented techniques can be adapted.

Our proof hinges on the fact that the SNA $\varphi^+$ can be
approximated by the iterates of the upper bounding line
$\T^1\times\{1\}$ of the phase space, whose geometry can be controlled
quite accurately. This observation has already been used in
\cite{jaeger:2007} and will be exploited further here. An outline of
the strategy is given in Section~\ref{Outline}. In
Section~\ref{Estimates} we derive the required estimates on the
approximating curves, which are used to compute the Hausdorff
dimension and the pointwise dimension in Section~\ref{Hausdorff}.
\bigskip

\noindent{\bf Acknowledgments.} This work was supported by an Emmy-Noether-Grant
of the German Research Council (DFG grant JA 1721/2-1) and is part of the
activities of the Scientific Network ``Skew product dynamics and multifractal
analysis'' (DFG grant OE 538/3-1). We thank Ren\'e Schilling for his thoughtful
comment leading to Proposition~\ref{p.pinchedset}.

\section{Preliminaries} \label{Preliminaries}

\subsection{Strange non-chaotic attractors} \label{SNA}

In the following, we provide some basics on SNA in pinched skew products by
sketching Keller's proof for the existence of SNA \cite{Keller1996}. More
precisely, according to \cite{grebogi/ott/pelikan/yorke:1984} the upper bounding
graph $\varphi^+$ is called an SNA if it is non-continuous and has a negative
Lyapunov exponent, and we will mainly explain how to obtain the non-continuity.

Let $I=[0,1]$ and $\T^D=\R^D/\Z^D$. A {\em quasiperiodically forced interval
map} is a skew product map of the form
\begin{equation*}
T \ : \ \T^D\times I\to\T^D \times I \quad , \quad (\theta,x)\mapsto
(\omega(\theta),T_\theta(x)) \ ,
\end{equation*}
where $\omega:\T^D\to\T^D,\ \theta\mapsto \theta+\rho\bmod 1$ an
irrational rotation. The maps $T_\theta:I\to I$ are called {\em fibre
  maps}. $T$ is {\em pinched} if there exists some
$\theta_*\in\T^D$ with $\#T_{\theta_*}(I) = 1$.

We denote by $\mathcal{T}$ the class of quasiperiodically forced
interval maps $T$ which share the following properties:
\begin{itemize}
\item[$(\mathcal{T}1)$] the fibre maps $T_\theta$ are monotonically increasing;
\item[$(\mathcal{T}2)$] the fibre maps $T_\theta$ are differentiable and
  $(\theta,x)\mapsto T_\theta'(x)$ is continuous on $\T^D\times I$;
\item[$(\mathcal{T}3)$] $T$ is pinched;
\item[$(\mathcal{T}4)$] $T_\theta(0)=0$ for all $\theta\in\T^D$.
\end{itemize}
Note that the last item means that the zero line $\T^D\times\{0\}$ is
$T$-invariant. It is easy to check that the maps $F_\kappa$ defined in
(\ref{e.gopy-map}) belong to $\mathcal{T}$.

An {\em invariant graph} is a measurable function $\varphi:\T^D\to I$
which satisfies (\ref{e.inv-graph}). If all fibre maps are
differentiable, the Lyapunov exponent of $\varphi$ is given by
$\lambda(\varphi):=\int_{\T^D}\log T_\theta'(\varphi(\theta)) \
d\theta$. The {\em upper bounding graph} $\varphi^+$ is given by
(\ref{e.upper-boundinggraph}). Equivalently, it can be defined by
\begin{equation*}
\varphi^+(\theta) \ = \ \nLim T^n_{\omega^{-n}(\theta)}(1),
\end{equation*}
where $T^n_\theta=T_{\omega^{n-1}(\theta)}\circ \ldots \circ T_\theta$. This
means that the {\em iterated upper bounding lines} 
\begin{equation}\label{e.upper-boundinglines}
\varphi_n(\theta) \ := \ T^n_{\omega^{-n}(\theta)}(1)
\end{equation}
converge pointwise and, by monotonicity of the fibre maps, in a decreasing way
to $\varphi^+$. This fact will be crucial for our later analysis. A first
consequence of this observation is that, under some mild conditions, the
Lyapunov exponent of $\varphi^+$ is always non-positive.

\begin{lem}[{\cite[Lemma 3.5]{jaeger:2003}}]\label{l.upper-lyap}
If $\theta\mapsto\log\left(\inf_{x\in I} T'_\theta(x)\right)$ is integrable,
then $\lambda(\varphi^+)\leq 0$.
\end{lem}

Now, turning back to the maps $F_\kappa$ defined in
(\ref{e.gopy-map}), the Lyapunov exponent of the zero line is easily
computed and one obtains
\begin{equation*}
\lambda(0) \ = \ \log \kappa - \log 2.
\end{equation*}
Consequently, when $\kappa>2$ this exponent is positive and therefore the upper
bounding graph cannot be the zero line. However, at the same time the pinching
condition together with the invariance of $\varphi^+$ imply that
$\varphi^+(\theta)=0$ for a dense set of $\theta\in\T^1$. Hence, $\varphi^+$
cannot be continuous.

Using the concavity of the fibre maps, it is further possible to show that
$\varphi^+$ is the only invariant graph of the system (\ref{e.gopy-map})
besides the zero line, that $\lambda(\varphi^+)$ is strictly negative and that
$\varphi^+$ attracts $\Leb_{\kreis\times I}$-a.e.\ initial condition
$(\theta,x)$, in the sense that
\begin{equation*}
  \nLim F_{\kappa,\theta}^n(x) - \varphi^+(\theta+n\rho\bmod 1) \ = \ 0. 
\end{equation*}
Finally, we note that to any invariant graph $\varphi$, an invariant measure
$\mu_\varphi$ can be associated by
\begin{equation*}
	\mu_{\varphi}(A) \ := \ \Leb_{\T^D}(\pi_1(A\cap\Phi))
\end{equation*}
for all Borel measurable sets $A\ssq \T^D\times I$, where
$\pi_1:\T^D\times I\to\T^D$ is the projection to the first coordinate.

\subsection{Dimensions} \label{Dimensions}

Let $X$ be a separable metric space.  The diameter of a subset $A\subseteq X$ is
denoted by $\diam(A)$.  For $\varepsilon>0$ a finite or countable collection
$\{A_i\}$ of subsets of $X$ is called an {\em $\varepsilon$-cover} of $A$ if
$\diam(A_i)\leq\varepsilon$ for each $i$ and $A\subseteq\bigcup_i A_i$.

\begin{definition}
	For $A\subseteq X$, $s\geq 0$ and $\varepsilon>0$ define
	\[
	\mathcal H_\varepsilon^s(A)\ :=\ \inf\left\{\left.\sum\limits_i
			(\diam(A_i))^s \ \right|\ \{A_i\}\text{ is an
			$\varepsilon$-cover of $A$}\right\}.
	\]
	Then
	\[
		\mathcal H^s(A)\ :=\ \lim\limits_{\varepsilon\to 0} \mathcal
			H_\varepsilon^s(A)
	\]
	is called the $s$-dimensional Hausdorff measure of $A$.  The Hausdorff
	dimension of $A$ is defined by
	\[
		D_H(A)\ :=\ \sup\{s\geq 0 \mid \mathcal H^s(A)=\infty\}.
	\]
\end{definition}

\begin{definition}
	The lower and upper box-counting dimension of a totally bounded subset
	$A\subseteq X$ are defined as
	\begin{align*}
		\underline D_B(A)\ :=\ \liminf\limits_{\varepsilon\to 0}
			\frac{\log N(A,\varepsilon)}{-\log\varepsilon},\\
		\overline D_B(A)\ :=\ \limsup\limits_{\varepsilon\to 0}
			\frac{\log N(A,\varepsilon)}{-\log\varepsilon},
	\end{align*}
	where $N(A,\varepsilon)$ is the smallest number of sets of diameter
	$\varepsilon$ needed to cover $A$.  If $\underline D_B(A)=\overline
	D_B(A)$, then their common value $D_B(A)$ is called the box-counting
	dimension (or capacity) of $A$.
\end{definition}

In general, we have $D_H(A)\leq D_B(A)$.
In the following, we will state some well known properties of the Hausdorff
measure and dimension that will be used later on.

\begin{lem}[\cite{Pesin1997}] \label{lemma_Hausdorff_dimension_Lipschitz_image}
	Let $X,Y$ be two separable metric spaces and let $g:A\subseteq X\to Y$
	be a Lipschitz continuous map with Lipschitz constant $K$.  Then
	$\mathcal H^s(g(A))\leq K^s\mathcal H^s(A)$ and $D_H(g(A))\leq D_H(A)$.
	Further, if $g$ is bi-Lipschitz continuous, then $D_H(g(A))=D_H(A)$.
\end{lem}

\begin{lem}[\cite{Pesin1997}]
	The Hausdorff dimension is countably stable, i.e.\ $D_H\left(\bigcup_i
	A_i\right)=\sup_i D_H(A_i)$ for any sequence of subsets $(A_i)_{i\in\set
	N}$ with $A_i\subseteq X$ .
\end{lem}

In contrast to the last lemma, we have that the upper box-counting dimension is only
finitely stable and that $D_B(A)=D_B\left(\overline A\right)$.

\begin{thm}[\cite{Howroyd1996}] \label{theorem_Hausdorff_dimension_product_sets}
	Let $X,Y$ be two separable metric spaces and consider the Cartesian
	product space $X\times Y$ equipped with the maximum metric.  Then for
	$A\subseteq X$ and $B\subseteq Y$ totally bounded we have 
	\[
		D_H(A\times B)\ \leq\ D_H(A)+\overline D_B(B).
	\] 
\end{thm}

\begin{lem} \label{lemma_Hausdorff_dimension_limsup_set}
	Let $A\subseteq X$ be a $\limsup$ set, meaning that there exists a sequence
	$(A_i)_{i\in\set N}$ of subsets of $X$ with
	\[
		A\ = \ \limsup\limits_{i\to\infty} A_i\ = \
		\bigcap\limits_{i=0}^\infty\bigcup\limits_{k=i+1}^\infty A_k.
	\]
	If $\sum_{i=1}^\infty\diam(A_i)^s<\infty$ for some $s>0$, then
	$\mathcal H^s(A)=0$ and $D_H(A)\leq s$.
\end{lem}
\proof
	Since $\sum_{i=1}^\infty\diam(A_i)^s<\infty$, we have 
	$\sum_{i=k}^\infty\diam(A_i)^s\to 0$ for $k\to\infty$.
	That means the diameter of the $A_i$'s goes to $0$ for $i\to\infty$.
	Therefore, $\{A_i: i\geq k\}$ is an $\varepsilon$-cover for $k$
	sufficiently large.
	This implies $\mathcal H_\varepsilon^s(A)\leq\sum_{i=k}^\infty\diam(A_i)^s
		\to 0$ for $k\to\infty$.
	Hence, $\mathcal H^s(A)=0$ and $D_H(A)\leq s$.
\qed\medskip

For $x\in X$ and $\varepsilon>0$ we denote by $B_\varepsilon(x)$
the open ball around $x$ with radius $\varepsilon>0$.

\begin{definition}
  Let $\mu$ be a finite Borel measure in $X$. For each point $x$ in
  the support of $\mu$ we define the lower and upper pointwise
  dimension of $\mu$ at $x$ as
	\begin{align*}
		\underline d_\mu(x)\ :=\ \liminf\limits_{\varepsilon\to 0}
			\frac{\log\mu(B_\varepsilon(x))}{\log\varepsilon},\\
		\overline d_\mu(x)\ :=\ \limsup\limits_{\varepsilon\to 0}
			\frac{\log\mu(B_\varepsilon(x))}{\log\varepsilon}.
	\end{align*}
	If $\underline d_\mu(x)=\overline d_\mu(x)$, then their common value 
	$d_\mu(x)$ is called the pointwise dimension of $\mu$ at $x$.
	We say that the measure $\mu$ is exact dimensional if the pointwise
	dimension exists and is constant almost everywhere, i.e.\
	\[
		\underline d_\mu(x)\ = \ \overline d_\mu(x)\ =:\ d_\mu,
	\]
	$\mu$-almost everywhere.
\end{definition}

\begin{definition}
	The lower and upper information dimension of $\mu$ are defined as
	\begin{align*}
		\underline D_1(\mu):=\liminf\limits_{\varepsilon\to 0}
			\frac{\int\log\mu(B_\varepsilon(x))d\mu(x)}{\log\varepsilon},\\
		\overline D_1(\mu):=\limsup\limits_{\varepsilon\to 0}
			\frac{\int\log\mu(B_\varepsilon(x))d\mu(x)}{\log\varepsilon}.
	\end{align*}
	If $\underline D_1(\mu)=\overline D_1(\mu)$, then their common value 
	$D_1(\mu)$ is called the information dimension of $\mu$.
\end{definition}

\begin{thm}[\cite{cutler1991dimensions,Zindulka2002}] 
     \label{theorem_relation_pointwise_information_dimension}
	Suppose $\overline D_B(X)<\infty$. We have
	\[
		\int\underline d_\mu(x)\ d\mu(x)\ \leq\ \underline D_1(\mu)
		\ \leq\ \overline D_1(\mu)\ \leq\ \int\overline d_\mu(x)\ d\mu(x).
	\]
	In particular, if $\mu$ is exact dimensional, then $D_1(\mu)=d_\mu$.
\end{thm}

Note that also several other dimensions of $\mu$ coincide if $\mu$ is
exact dimensional
\cite{young1982dimension,Zindulka2002,ledrappier1985metric,ledrappier1985metric2}.

\subsection{Rectifiable sets and measures} \label{rectifiability}

Here, we mainly follow \cite{AmbrosioKirchheim2000}.

\begin{definition}
	For $D\in\set N$ a Borel set $A\subseteq X$ is called countably
	$D$-rectifiable if there exists a sequence of Lipschitz continuous functions
	$(g_i)_{i\in\set N}$ with $g_i:A_i\subseteq\set R^D\to X$ such that 
	$\mathcal H^D(A\backslash\bigcup_i g_i(A_i))=0$.
	A finite Borel measure $\mu$ is called $D$-rectifiable if
	$\mu=\Theta\left.\mathcal H^D\right|_A$ for some countably $D$-rectifiable
	set $A$ and some Borel measurable density $\Theta:A\to[0,\infty)$.
\end{definition}

Note that, by the Radon-Nikodym theorem, $\mu$ is $D$-rectifiable if and only
if $\mu$ is absolutely continuous with respect to $\left.\mathcal H^D\right|_A$
with $A$ some countably $D$-rectifiable set.

\begin{thm}[{\cite[Theorem 5.4]{AmbrosioKirchheim2000}}]
	For a $D$-rectifiable measure $\mu=\Theta\left.\mathcal H^D\right|_A$ we have
	\[
		\Theta(x) \ = \ \lim\limits_{\varepsilon\to 0}
		\frac{\mu(B_\varepsilon(x))}{V_D\varepsilon^D},
	\]
	for $\mathcal H^D$-a.e.\ $x\in A$, where $V_D$ is the volume of the
	$D$-dimensional unit ball.
	The right hand side of this equation is called $D$-density of $\mu$.
\end{thm}

This theorem implies in particular that the $D$-density exists and is positive
$\mu$-almost everywhere for a $D$-rectifiable measure $\mu$ and this gives directly

\begin{cor}\label{dimensions_rectifiable_measure}
	A $D$-rectifiable measure $\mu$ is exact dimensional with $d_\mu=D_1(\mu)=D$.
\end{cor}

\section{Outline of the strategy} \label{Outline}

As we have mentioned in the introduction, our main goal is to analyze the
structure of the upper bounding graphs $\varphi^+$ when they are different from
the zero line, and in particular we want to determine the dimensions of these
graphs and of their associated invariant measures.  However, the argument
for the non-continuity of the invariant graphs sketched in Section~\ref{SNA} is
a `soft' one and does not yield any quantitative information about the structure
of the invariant graphs. Hence, it is not clear how such an analysis can be
carried out.

\begin{figure}[h] 
\begin{center}
\epsfig{file=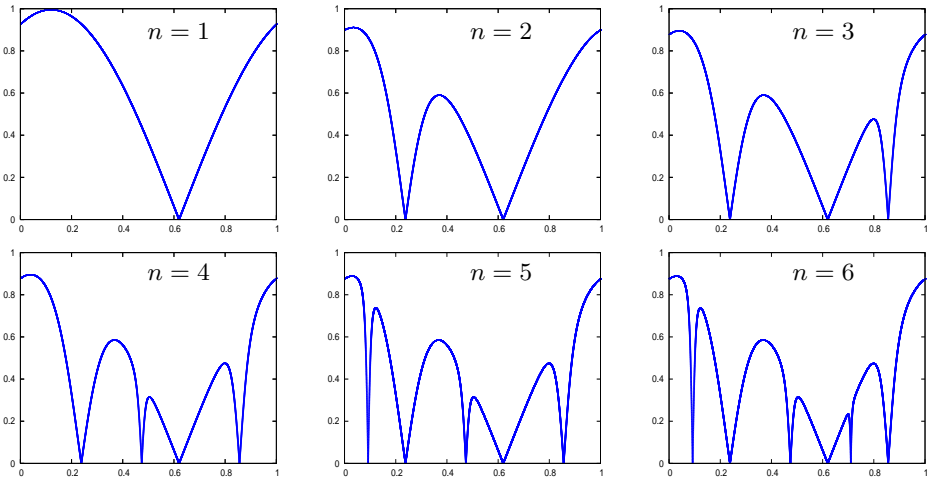, width=0.9\linewidth}
\end{center} {\caption{ \small The graphs of the first six iterated
    upper bounding lines of (\ref{e.gopy-map}) with $\kappa=3$ and
    $\rho$ the golden mean. \label{f.2}}}
\end{figure}

However, as mentioned above the upper bounding graph $\varphi^+$ can
be approximated by the iterated upper bounding lines $\varphi_n$
defined in (\ref{e.upper-boundinglines}). It turns out that the
geometry of the lines $\varphi_n$ can be controlled well, and this is
the starting point of our investigation.  Figure~\ref{f.2} shows the
first six iterates $\varphi_1 \ld \varphi_6$. A clear pattern can be
observed.  Apparently, when going from $\varphi_{n-1}$ to
$\varphi_{n}$, the only significant change is the appearance of a new
`peak' in a small ball $I_n$ around the $n$-th iterate
$\tau_n=\omega^n(\theta_*)$ of the pinching point $\theta_*$.  Outside
of $I_n$, the graphs seem to remain unchanged. Further, since every
new peak is the image of the previous one and due to the expansion
around the 0-line, the peaks become steeper and sharper in every step.
As a consequence, the radius of the balls $I_n$ decreases
exponentially.

Of course, this is a very rough picture, which can only hold in an approximate
sense. Due to the strict monotonicity of the fibre maps for all $\theta\neq
\theta_*$, the sequence $\varphi_n$ is strictly decreasing everywhere except on
the countable set $\{\tau_n\mid n\in\N\}$, so the graphs have to change at least
a little bit outside of $I_n$. However, let us assume for the moment that the
above description was true and $\varphi_{n-1}(\theta)-\varphi_n(\theta)=0$ for
all $\theta\notin I_n$. In this case, the graph $\varphi^+$ is already
determined on $\T^D\smin\bigcup_{k=n}^\infty I_k=:\Lambda_n$ after $n$ steps and
equals $\varphi_{n|\Lambda_n}$ on this set. However, as a finite iterate of
$\T^D\times\{1\}$, the function $\varphi_n$ is Lipschitz continuous and
therefore its graph $\Phi_{n|\Lambda_n}=\{(\theta,\varphi_n(\theta)) \mid
\theta\in\Lambda_n\}$ has Hausdorff dimension $D$. Due to the exponential
decrease of the radius of the $I_n$, the set
$\Omega_\infty=\T^D\smin\bigcup_{n\in\N} \Lambda_n$ is a $\limsup$ set and has
Hausdorff dimension zero by Lemma~\ref{lemma_Hausdorff_dimension_limsup_set}. It
follows that $\Phi^+$ is contained in the countable union $\bigcup_{n\in\N}
\Phi_{n|\Lambda_n}\cup (\Omega_\infty\times [0,1])$ of at most D-dimensional
sets. By countable stability, this implies that the Hausdorff dimension of
$\Phi^+$ is $D$. For the pointwise dimension, a similar argument could be given
but we will directly conclude from the arguments sketched above that
$\mu_{\varphi^+}$ is $D$-rectifiable.

The remainder of this article is devoted to showing that these heuristics can be
converted into a rigorous proof, despite the fact that `nothing changes outside
of $I_n$' has to be replaced by `almost nothing changes outside of $I_n$'.

\section{Estimates on the iterated upper bounding lines}\label{Estimates}

The purpose of this section is to obtain a good control on the
behaviour and shape of the iterated upper bounding lines. In order to
derive the required estimates, we have to impose a number of
assumptions on the geometry of our systems. The hypotheses are
formulated in terms of ${\cal C}^1$-estimates, and it is easy to check
that they are fulfilled by (\ref{e.gopy-map}) whenever $\kappa$ is
large enough (see Lemma~\ref{l.estimate_check} for details).
\smallskip

Let $T\in\mathcal{T}$. Suppose there exist $\alpha>2,\ \gamma>0$ and
$L_0\in (0,1)$ such that for all $\theta\in\T^D$
\begin{equation} \label{Lipschitz-fibres}
	\abs{T_\theta(x)-T_\theta(y)}\ \leq \ \alpha\abs{x-y},
\end{equation}
for all $x,y\in[0,1]$, and
\begin{equation} \label{Lipschitz-contraction-fibers}
	\abs{T_\theta(x)-T_\theta(y)}\ \leq \ \alpha^{-\gamma}\abs{x-y},
\end{equation}
for all $x,y\in[L_0,1]$.
Further, we assume there exists $\beta>0$ such that for all $x\in[0,1]$
\begin{equation} \label{Lipschitz_condition_base}
	\abs{T_\theta(x)-T_{\theta'}(x)}\ \leq\ \beta d(\theta,\theta').
\end{equation}
When $T$ is differentiable in $\theta$, we may for example take
$\beta=\sup_{(\theta,x)}\Abs{\partial_\theta T_\theta(x)}$. As above, we let
$\tau_n:=\omega^n(\theta_*)$. We suppose the rotation vector $\rho\in\R^D$ is
Diophantine, meaning that there exist constants $c>0$ and $d>1$ such that 
\begin{equation} \label{Diophantine_condition}
	d(\tau_n,\theta_*)\ \geq \ c\cdot n^{-d},
\end{equation}
for all $n\in\N$. In addition, we assume there are $m\in\set N$, $a>1$ and
$0<b<1$ with
\begin{align}
	m& >\ 22\left(1+\frac{1}{\gamma}\right),\label{condition_constant_m}\\ a&
	\ \geq\ (m+1)^d, \label{condition_constant_a}\\ b& \ \leq\
	c,\label{condition_constant_b} \\ d(\tau_n,\theta_*)&\ > \ b\quad
	\textnormal{ for all }n\in\{1,\dots,m-1\}
	\label{pinched_orbit_condition_first_iterates}
\end{align}
such that 
\begin{align} \label{definition_reference_system}
	T_\theta(x)&\ \geq\
	\min\{L_0,ax\}\cdot\min\left\{1,\frac{2}{b}d(\theta, \theta_*)\right\},
\end{align}
for all $(\theta,x)\in\T^D\times[0,1]$. We now let 
\begin{equation}
  \label{eq:1}
  \mathcal{T}^* \ := \ \left\{ T\in\mathcal{T} \mid T \textrm{ satisfies
  (\ref{Lipschitz-fibres})--(\ref{definition_reference_system})}\right\} \ ,
\end{equation}
where ``satisfies
(\ref{Lipschitz-fibres})--(\ref{definition_reference_system})''
should be understood in the sense of ``there exist constants
$\alpha$, $\gamma$, $L_0$, $\beta$, $c$, $d$, $m$, $a$ and $b$ such that
(\ref{Lipschitz-fibres})--(\ref{definition_reference_system}) are satisfied''.
\begin{example} \label{example_higher_dimensional_case} The following
  map is a simple extension of \eqref{e.gopy-map} with a
  higher-dimensional rotation on the base. 
	\begin{equation} \label{e.F_kappa}
          F_\kappa \ : \ \T^D \times [0,1] \to \T^D \times [0,1] \quad,\quad
          F_\kappa(\theta,x) \ = \ \left(\theta+\rho\bmod 1,
            \tanh\left(\kappa x\right)\cdot \frac{1}{D} \cdot\sum\limits_{i=1}^D
            \sin(\pi\theta_i)\right)
	\end{equation}
	Here $\theta=(\theta_1,\dots,\theta_D)$. 
\end{example}
As we show now, $F_\kappa$ satisfies \eqref{Lipschitz-fibres} --
\eqref{definition_reference_system} for all sufficiently large
$\kappa$.
\begin{lem}\label{l.estimate_check}
  Let $\rho$ satisfy the Diophantine condition
  (\ref{Diophantine_condition}) with constants $c,d$. Then there exist
  constants $D_0=D_0(c,d)$ and $\kappa_0=\kappa_0(c,d,D)$ such that
  \begin{itemize}
  \item for all $\kappa\geq \kappa_0$ the map $F_\kappa$ belongs to
    $\mathcal{T}^*$;
  \item if $D\geq D_0$, then the constants $\alpha$, $m$ and $a$ can be
    chosen such that
    \begin{equation}
      \label{eq:2}
      D \ > \ m^2 \log(\alpha/a) \ .
    \end{equation}
  \end{itemize}
 \end{lem}
 The additional condition (\ref{eq:2}) will be used to show that for
 sufficiently large $D$ the $D$-dimensional Hausdorff measure of the
 upper bounding graph $\varphi^+$ of $F_\kappa$ is finite, see
 Proposition~\ref{proposition_Hausdorff_measure_finite}.  

 \proof We let $\alpha=\kappa$, $\gamma=\halb$,
 $L_0=\frac{\log\kappa}{\kappa}$, $\beta=\pi$, $m=67$,
 $b=\halb\min_{n=1}^{m-1} cn^{-d}$ and
 $a=\frac{2b\kappa}{D(e+1/e)^2}$. Then we choose $D_0=D_0(c,d)$ such
 that for all $D\geq D_0$
 \begin{equation}
   \label{e.D0-cond}
   D \ > \ m^2\log\left(\frac{D(e+1/e)^2}{2b}\right)
 \end{equation}
 and $\kappa_0=\kappa_0(c,d,D)$ such that for all $\kappa\geq
 \kappa_0$
 \begin{eqnarray}
   \label{eq:4}
   \kappa & \geq & 16 \ , \\ \label{e.4b}
   \frac{2b\kappa}{D(e+1/e)^2} & \geq & (m+1)^d \ , \\
   \frac{\log\kappa}{\kappa} & \leq & \frac{b\tanh(1)}{2D} \ .\label{e.4c}
 \end{eqnarray}
We have 
\begin{equation}
  \label{eq:3}
  \left[\tanh(\kappa x)\right]' \ = \ \frac{4\kappa}{(e^{\kappa x}+e^{-\kappa x})^2} \ \leq \ \kappa
\end{equation}
for all $x\geq 0$ and 
\begin{equation}
  \label{eq:5}
  0 \ \leq \ \frac{1}{D}\sum_{i=1}^D \sin(\pi\theta_i) \ \leq \ 1
\end{equation}
for all $\theta\in\T^D$. Hence, (\ref{Lipschitz-fibres}) holds and
since
\begin{equation}
  \label{eq:6}
  F'_{\kappa,\theta}(x)\ \leq\  F'_{\kappa,\theta}(L_0) \ \leq \ 
  \frac{4\kappa}{(\kappa+1/\kappa)^2} 
  \ \leq \ \frac{4}{\kappa} \ \leq \ \kappa^{-1/2}
\end{equation}
for all $x\geq L_0$, the same is true for
(\ref{Lipschitz-contraction-fibers}). (\ref{Lipschitz_condition_base})
and (\ref{condition_constant_m}) are easy to check, and
(\ref{Diophantine_condition}) holds by assumption.
(\ref{condition_constant_a}) follows from (\ref{e.4b}), whereas
(\ref{condition_constant_b}) and
(\ref{pinched_orbit_condition_first_iterates}) are obvious from the
choice of $b$ and (\ref{Diophantine_condition}).  In order to verify
(\ref{definition_reference_system}), note that $[\tanh(\kappa
x)]'_{|x=1/\kappa} = \frac{4\kappa}{(e+1/e)^2}$, such that by
concavity and monotonicity
\begin{equation}
 \tanh(\kappa x) \ \geq \ \left\{
   \begin{array}{ll} \frac{4\kappa}{(e+1/e)^2}\cdot x & \textrm{if } x\leq 1/\kappa \\ \\
      \tanh(1) & \textrm{if } x>1/\kappa 
        \end{array} \right. \ .
\end{equation}
Using (\ref{e.4c}) and the fact that $\sum_{i=1}^D \sin(\pi\theta_i) \geq
d(\theta,\theta_*)$, where $\theta_*=0$, we obtain
\begin{eqnarray*}
  \label{eq:7}
  F_{\kappa,\theta}(x) & \geq & \min\left\{ \tanh(1),
    \frac{4\kappa}{(e+1/e)^2} x\right\} \cdot \frac{1}{D}
   d(\theta,\theta_*) \\
  & \geq & \min\left\{ \frac{b\tanh(1)}{2D}, \frac{2b\kappa}{D(e+1/e)^2} x\right\} 
 \cdot \frac{2}{b}d(\theta,\theta_*) \\
 & \geq & \min\{L_0,ax\} \cdot \min\left\{1,\frac{2}{b}d(\theta,\theta_*)\right\} \ 
 \end{eqnarray*}
 as required. Finally, since $\alpha/a = \frac{D(e+1/e)^2}{2b}$
 condition (\ref{eq:2}) follows from (\ref{e.D0-cond}). Note that
 since $b$ and $m$ are constants only depending on $c$ and $d$, the
 same is true for the condition~(\ref{e.D0-cond}) on $D_0$.
 \qed\medskip

\begin{rem}
\label{r.0-lyap}
Given $T\in\mathcal{T}^*$, note that
(\ref{definition_reference_system}) implies 
\[
	\lambda(0)\ \geq \
	\log \frac{2a}{b}+\int_{\T^D}\log d(\theta,\theta_*) \ d\theta\ \geq\
	\log\frac{2a}{b}-\log 2-1\ .
\]
Since $a\geq 23$ by (\ref{condition_constant_a}), this yields $\lambda(0)>0$
and hence $\varphi^+(\theta)>0$ for $\Leb_{\T^D}$-almost every $\theta$.
\end{rem}
\medskip

In order to formulate the main results of this section, let $j\in\R$ and 
\[
	r_j \ := \ \frac{b}{2}a^{-\frac{j-1}{m}}\ .
\]

\begin{prop}\label{properties_approx_graphs}
  Let $T\in\mathcal{T}^*$. Given $q\in\N$, the following hold.
\begin{enumerate}
 \item[(i)] $|\varphi_{n}(\theta)- \varphi_{n}(\theta')|\leq\beta\alpha^n
	d(\theta,\theta')$ for all $n\in\N$ and $\theta,\theta'\in\T^D$.
      \item[(ii)] There exists $\lambda>0$ such that if $n\geq mq+1$ and
        $\theta\notin\bigcup_{j=q}^n B_{r_{j}}(\tau_j)$, then
        $|\varphi_{n}(\theta)-\varphi_{n-1}(\theta)|\leq
        \alpha^{-\lambda(n-1)}$.
 \item[(iii)] There exists $K>0$ such that if
        $\theta,\theta'\notin\bigcup_{j=q}^n B_{r_j}(\tau_j)$, then
        $|\varphi_{n}(\theta)- \varphi_{n}(\theta')|\leq
        K\alpha^{mq}d(\theta,\theta')$ for all $n\in\set N$. 
\end{enumerate}
\end{prop}
For the proof, we need two preliminary statements. The first is a simple
observation.

\begin{lem}\label{lemma_two_assertions}
  Suppose (\ref{Diophantine_condition}) holds and let $n,i\in\set N_0$
  and $n>0$. If $d(\tau_n,\theta_*)\leq b\cdot a^{-i}$, then $n\geq
  a^{i/d}$.
\end{lem}
\proof \eqref{Diophantine_condition} implies $c\cdot n^{-d}\leq
		b\cdot a^{-i}$, and using \eqref{condition_constant_b} we get
		$n^{-d}\leq a^{-i}$.
\qed\medskip

The second statement we need for the proof of 
Proposition~\ref{properties_approx_graphs} is an upper bound on the proportion
of time the backwards orbit of a point $(\theta,\varphi_n(\theta))\in\Phi_n$
spends outside of the contracting region $\T^D\times[L_0,1]$.
Given $\theta\in\T^D$ and $n\in\N$, let $\theta_k:=\omega^{k-n}(\theta)$ and
$x_k:=\varphi_k(\theta_k)$ for $0\leq k\leq n$.
Note that thus $x_k=T^k_{\theta_0}(1)$ and 
$T_{\theta_k}^{n-k}(x_k)=\varphi_n(\theta)$.
Let
\begin{eqnarray*}
 s^n_k(\theta) & := & \#\{k\leq j < n \mid x_j < L_0\} \quad \textrm{and}\\
 s^n_k(\theta,\theta') & := & \#\{k \leq j < n \mid \min\{x_j,x_j'\} < L_0\} \
\end{eqnarray*}
and note that $s^n_k(\theta,\theta')\leq s^n_k(\theta)+s^n_k(\theta')$. We set
$s^n_n(\theta):=0$ and $s^n_n(\theta,\theta'):=0$.

\begin{lem}\label{main_property_of_s_n} Let $T\in\mathcal{T}^*$ and 
  $q,n\in\N$ with $n\geq mq+1$. Suppose that
  $\theta\notin\bigcup_{j=q}^n B_{r_j}(\tau_j)$.  Then for all $t\geq
  mq$ we have
	\[
		s^n_{n-t}(\theta)\leq \frac{11t}{m} \ .
	\]
\end{lem}
\proof
We divide $A=\{1\leq k< n-q \mid x_k< L_0\}$ into blocks $B=\{l+1,\dots,p\}$
with $0\leq l<p<n-q$ and the properties
\alphlist
	\item $x_l\geq L_0/a$;
	\item $x_k<L_0/a$ for all $k\in\{l+1,\dots,p-1\}$;
		\item $x_p<L_0$;
	\item either $x_p\geq L_0/a$ or $x_{p+1}\geq L_0$ or $p+1=n-q$.
\listend 

Note that these blocks cover the whole set $A$, and they are uniquely determined
by the above requirements. Since we always start a new block when the
`threshold' $L_0/a$ is reached, we may have $p=l'$ for two adjacent blocks
$B=\{l+1\ld p\}$ and $B'=\{l'+1\ld p'\}$. 

Now, we first consider a single block $B=\{l+1,\dots,p\}$.  We have $\theta_l\in
B_{b/2}(\theta_*)$, because otherwise $x_{l+1}\geq L_0$ according to
\eqref{definition_reference_system} and (a).  Since $x_{k+1}=T_{\theta_k}(x_k)$,
we can use \eqref{definition_reference_system} and (b) to obtain that for any
$k\in\{l+1,\dots,p-1\}$
\[
     x_{k+1}\geq ax_k\min\left\{1,\frac{2}{b}d(\theta_k,\theta_*)\right\}.
\]
Therefore, using (c),(a) and \eqref{definition_reference_system} again, we see that
\begin{align}
	\label{inequality_product_distances}
	1>\frac{x_p}{L_0}\geq a^{p-l-1}\prod\limits_{k=l}^{p-1}
			\min\left\{1,\frac{2}{b}d(\theta_k,\theta_*)\right\}\ .
\end{align}
Now, note that 
\begin{eqnarray*}
  \lefteqn{\sum\limits_{k=l}^{p-1}\log\min\left\{1,\frac{2}{b}d(\theta_k,
      \theta_*)\right\} \ \geq} \\ & \geq & -\ (\log a)\cdot 
 \sum_{i=1}^\infty i\cdot\#\left\{l\leq k<p\left|\ \frac{b}{2}a^{-i} \leq
      d(\theta_k,\theta_*)<\frac{b}{2}a^{-i+1}\right.\right\}\ . \nonumber
\end{eqnarray*}
Therefore, we can deduce from \eqref{inequality_product_distances} that 
\begin{eqnarray} \label{upper_bound_length_block}
	p-l & \leq & \sum_{i=1}^\infty i\cdot\#\left\{l\leq k<p\left|\
	\frac{b}{2}a^{-i} \leq
	d(\theta_k,\theta_*)<\frac{b}{2}a^{-i+1}\right.\right\} \\ & = &
         \sum_{i=1}^\infty  \#\left\{l\leq k<p\left|\
	d(\theta_k,\theta_*)<\frac{b}{2}a^{-i+1}\right.\right\} \ . \nonumber
\end{eqnarray}
We turn to the estimate on $A \cap [n-t,n-q)$ (note that $n-t<n-q$).
It may happen that $n-t$ is contained in a middle of a block $B$.
In this case, we need two auxiliary statements to estimate the length of this
first block intersecting $[n-t,n-q)$.
Let $j\in\set N$ be such that $(m-3)(j-1)<t\leq (m-3)j$.
\begin{claim}
\label{claim1}
If $j'\geq 1$ and $d(\theta_k,\theta_*)\geq ba^{-j'}/2$ for all $k=l\ld p-1$, then 
$p-l\leq\frac{j'}{1-2/m} \leq 3j'$.
\end{claim}
\proof Due to (\ref{pinched_orbit_condition_first_iterates}), two consecutive
visits in $B_{b/2}(\theta_*)$ are at least $m$ times apart, whereas two
consecutive visits in $B_{b a^{-i}/2}(\theta_*)$ are at least $a^{i/d}$ times
apart by Lemma~\ref{lemma_two_assertions}. Hence, we obtain from
(\ref{upper_bound_length_block}) that 
\[
p-l \ \leq \ \frac{p-l}{m} +1 + \sum_{i=2}^{j'} \left(\frac{p-l}{a^{(i-1)/d}}
  +1\right) \ \stackrel{(\ref{condition_constant_a})}{\leq} \ \frac{2(p-l)}{m}
+j'.
\]
\roundqed\medskip

\begin{claim}\label{claim2}
Suppose the block $B=\{l+1\ld p\}$ intersects $[n-t,n-q)$ and $t\leq (m-3)j$.
Then $d(\theta_k,\theta_*)\geq ba^{-j+1}/2$ for all $k\in B$.
\end{claim}
\proof Suppose for a contradiction that there exist $j'\geq j$ and $k'\in B$
with $d(\theta_{k'},\theta_*) < ba^{-j'+1}/2$.
If $j'$ is chosen maximal, such that $d(\theta_k,\theta_*)\geq ba^{-j'}/2$ for
all $k\in B$, then Claim~\ref{claim1} implies that $\#B\leq 3j'$.
However, since $\theta\notin\bigcup_{k=q}^n B_{r_k}(\tau_k)$ we have
$d(\theta_k,\theta_*)\geq r_{n-k}$ for all $k\in\{0,\dots,n-q\}$ and this
implies $ba^{-j'}/2\geq r_{n-k'}$, i.e.\ $k'<n-mj'$.
Therefore, $n-t \leq \max B \leq k'+3j' < n-(m-3)j'$, contradicting the
assumption on $t$.  \roundqed\medskip

We can now complete the proof of the lemma. 
For all blocks $B$ intersecting $[n-t,n-q)$, Claim~\ref{claim2} implies
$d(\theta_k,\theta_*)\geq ba^{-j+1}/2$ for all $k\in B$, such that
$\# B\leq 3j$ by Claim~\ref{claim1}. Hence, by the same counting
argument as in the proof of Claim~\ref{claim1} and summing up over all blocks,
we obtain the following estimate from (\ref{upper_bound_length_block})
\begin{eqnarray*}
  s^n_{n-t}(\theta) & \leq & q + \#(A\cap [n-t,n-q))  \\
  & \leq & q + 3j + \frac{t}{m} + 1 + \sum_{i=2}^{j-1}\frac{t}{a^{(i-1)/d}}+1  \\
  & \leq & q + 4j + \frac{2t}{m}\ \ \stackrel{(\ref{condition_constant_m})}{\leq} 
  \ \ \frac{11t}{m}
\end{eqnarray*}
(recall that $t\geq mq$).
\qed\medskip

This allows to turn to the
\proof[\bf Proof of Proposition~\ref{properties_approx_graphs}] (i) For all
$\theta,\theta'\in\T^D$, we have 
\begin{align}\label{Lipschitz_property_induction_basis}
        \abs{\varphi_1(\theta)-\varphi_1(\theta')} \ = \
        \abs{T_{\omega^{-1}(\theta)}(1)-T_{\omega^{-1}(\theta')}(1)} \
        \stackrel{~(\ref{Lipschitz_condition_base})}{\leq} \ \beta
        d(\omega^{-1}(\theta),\omega^{-1}(\theta')) \ = \ \beta 
        d(\theta,\theta')
\end{align}
and
\begin{align}\label{Lipschitz_property_induction_step}
	\abs{\varphi_{n+1}(\theta)-\varphi_{n+1}(\theta')}
	\leq\ \abs{T_{\theta_n}(x_n)-T_{\theta_n}(x_n')
            \vphantom{T_{\theta_n'}(x_n')}}
	+\abs{T_{\theta_n}(x_n')-T_{\theta_n'}(x_n')}.
\end{align}
We claim that for all $\theta,\theta'\in\T^D$ 
\begin{equation} \label{e.simple_Lipschitz_induction}
    \abs{\varphi_n(\theta)-\varphi_n(\theta')} \ \leq\
    \beta(\alpha^n-1) d(\theta,\theta').
\end{equation}
For the proof of this assertion, we proceed by induction.
(\ref{e.simple_Lipschitz_induction}) holds for $n=1$ because of
\eqref{Lipschitz_property_induction_basis} and the fact that $\alpha>2$.
Moreover,
\begin{eqnarray*}
        \lefteqn{\abs{\varphi_{n+1}(\theta)-\varphi_{n+1}(\theta')} \
        \leq} \\ & \stackrel{(\ref{Lipschitz_property_induction_step})}{\leq} &
        \abs{T_{\omega^{-1}(\theta)}(\varphi_n(\omega^{-1}(\theta)))-
        T_{\omega^{-1}(\theta)}(\varphi_n(\omega^{-1}(\theta')))} \\ & & + \
        \abs{T_{\omega^{-1}(\theta)}(\varphi_n(\omega^{-1}(\theta')))-
        T_{\omega^{-1}(\theta')}(\varphi_n(\omega^{-1}(\theta'))) }
        \\ & \stackrel{~(\ref{Lipschitz-fibres}),(\ref{Lipschitz_condition_base})}{\leq}
        & \alpha  |\varphi_n(\theta')-\varphi_n(\theta)| \ +\  \beta  d(\theta,\theta') 
        \\ & \stackrel{(\ref{e.simple_Lipschitz_induction})}{\leq} & 
        \left(\alpha\beta(\alpha^n-1)+\beta\right)d(\theta,\theta')\ \
        \leq\ \ \beta(\alpha^{n+1}-1) d(\theta,\theta'),
\end{eqnarray*}
which proves (\ref{e.simple_Lipschitz_induction}) for $n+1$. \medskip

\noindent
(ii) We fix $n\in\N$ and $\theta\in\T^D$.
Let $\theta_k$ and $x_k$ be defined as above.
If $\varphi_{k-1}(\theta_k)-\varphi_k(\theta_k)=0$ for some $k\in\{1,\dots,n\}$,
then $\varphi_{n-1}(\theta_n)-\varphi_n(\theta_n)=0$.
Thus, we may assume that the distance is greater than $0$ for all $k$.
In this case, we have
\begin{eqnarray*}
  \lefteqn{\varphi_{n-1}(\theta)-\varphi_n(\theta)\ =  \ 
    (\varphi_{0}(\theta_1)-\varphi_1(\theta_1))\cdot \prod\limits_{k=1}^{n-1}
    \frac{\varphi_{k}(\theta_{k+1})-\varphi_{k+1}(\theta_{k+1})}
    {\varphi_{k-1}(\theta_k)-\varphi_k(\theta_k)}} \\ &\leq &
  \prod\limits_{k=1}^{n-1}
  \frac{T_{\theta_k}(\varphi_{k-1}(\theta_{k}))
    -T_{\theta_k}(\varphi_{k}(\theta_{k}))}
  {\varphi_{k-1}(\theta_k)-\varphi_k(\theta_k)} \ \leq \  
   \alpha^{s^n_1(\theta)-\gamma(n-1-s^n_1(\theta))},
\end{eqnarray*}
where we used \eqref{Lipschitz-fibres} and
\eqref{Lipschitz-contraction-fibers}.  Since
$\theta\notin\bigcup_{j=q}^n B_{r_j}(\tau_j)$, we can use Lemma
\ref{main_property_of_s_n} with $t=n-1$ to obtain
$|\varphi_n(\theta)-\varphi_{n-1}(\theta)|\leq
\alpha^{-\lambda(n-1)}$ where
\[
	\lambda:=\gamma-\frac{11}{m}(1+\gamma)
	\stackrel{(\ref{condition_constant_m})}{>}0.
\]

\noindent
(iii) We proceed by induction to show that for all $\theta,\theta'\in\T^D$ and
$n\in\N$ we have
\begin{equation} \label{e.Lipschitz_induction}
        \abs{\varphi_n(\theta)-\varphi_n(\theta')} \ \leq \ \beta
        \left(\sum\limits_{k=0}^{n-1}\alpha^{(1+\gamma)s^n_{n-k}(\theta,\theta')-\gamma
            k}\right) d(\theta,\theta').
\end{equation}
For $n=1$ this is true because of
\eqref{Lipschitz_property_induction_basis}.  Further, since
\begin{equation}\label{e.sn_sum}s^{n+1}_n(\theta,\theta') +
        s^n_{n-k}(\omega^{-1}(\theta),\omega^{-1}(\theta'))\ =\
        s^{n+1}_{n-k}(\theta,\theta'), 
\end{equation}
we have
\begin{eqnarray*}
  \lefteqn{\abs{\varphi_{n+1}(\theta)-\varphi_{n+1}(\theta')}
    \ \leq}
  \\ &  \stackrel{\eqref{Lipschitz_property_induction_step},
    \eqref{Lipschitz-fibres}-\eqref{Lipschitz_condition_base}}{\leq} 
    & \ \alpha^{(1+\gamma)s^{n+1}_n(\theta,\theta')-\gamma}
  \abs{\varphi_n(\omega^{-1}(\theta))-\varphi_n(\omega^{-1}(\theta'))} 
  +\beta d(\omega^{-1}(\theta),\omega^{-1}(\theta'))  \\
  &\stackrel{~(\ref{e.Lipschitz_induction}),(\ref{e.sn_sum})}{\leq}
   &\beta\left(\sum\limits_{k=0}^n\alpha^{(1+\gamma)
      s^{n+1}_{n+1-k}(\theta,\theta')-\gamma k} \right) d(\theta,\theta').
\end{eqnarray*}
This completes the induction step, such that (\ref{e.Lipschitz_induction}) holds
for all $n\in\N$.  \medskip

Now, when $\theta,\theta'\notin \bigcup_{j=q}^n B_{r_j}(\tau_j)$ and $k\geq mq$,
then $s^n_{n-k}(\theta,\theta')\leq \frac{22k}{m}$ by
Lemma~\ref{main_property_of_s_n}. Consequently, (\ref{e.Lipschitz_induction})
yields that 
\begin{eqnarray*}
  \abs{\varphi_n(\theta)-\varphi_n(\theta')}
  &\leq&\beta\left(\sum\limits_{k=0}^{mq-1}\alpha^k
    +\sum\limits_{k=mq}^{n-1}\alpha^{(1+\gamma)s^n_{n-k}(\theta,\theta')
	-\gamma k}\right)d(\theta,\theta')\\
  &\leq&\beta\left(\alpha^{mq}
    +\sum\limits_{k=mq}^{n-1}\alpha^{-\left(\gamma-\frac{22}{m}(1+\gamma)\right)k}
  \right)d(\theta,\theta').
\end{eqnarray*}
Because of \eqref{condition_constant_m}, we have
$\gamma-\frac{22}{m}(1+\gamma)>0$, and this implies
$|\varphi_{n}(\theta)-\varphi_{n}(\theta')|\leq K\alpha^{mq}d(\theta,\theta')$
with
\[
	K \ := \ \beta\left(1+\frac{\alpha^{-mq}}{1-
	\alpha^{-(\gamma-\frac{22}{m}(1+\gamma))}}\right). \qedhere
\]

\section{Dimensions of the upper bounding graph and the associated
  physical measure}
  \label{Hausdorff}

  For $T\in\mathcal{T}^*$, we can now calculate the Hausdorff
  dimension of the upper bounding graph $\varphi^+$, or more precisely
  of the corresponding point set $\Phi^+$.  We will also be able to
  draw some conclusions regarding the Hausdorff measure of $\Phi^+$.
  To that end, we will partition $\varphi^+$ into countably many
  subgraphs.  First, keeping the notation from the last section we
  define a partition of $\T^D$ by subsets $\Omega_j\subset\T^D$ with
  $j\in\set N_0\cup\{\infty\}$ as
\begin{eqnarray}
  \Omega_0 & := & \T^D\backslash\bigcup\limits_{k=j_0}^\infty B_{r_k}(\tau_k) \ ,\\
  \Omega_j & := &  B_{r_{j+j_0-1}}(\tau_{j+j_0-1})\backslash
  \bigcup\limits_{k=j+j_0}^\infty B_{r_k}(\tau_k) \ , \\
  \Omega_\infty & := & \bigcap\limits_{i=0}^\infty\bigcup\limits_{k=i+1}^\infty
  B_{r_k}(\tau_k) \ , \label{e.omega_infty}
\end{eqnarray}
where we choose $j_0\in\set N$ large enough to ensure $\Leb_{\T^D}(\Omega_j)>0$
for all $j\in\set N_0$.
This works for $j=0$ because $\sum_{k=1}^\infty\Leb_{\T^D}(B_{r_k}(\tau_k))<\infty$
and for $j\in\set N$ because for all $j'>j$ with
$B_{r_j}(\tau_j)\cap B_{r_{j'}}(\tau_{j'})\neq\emptyset$ the Diophantine condition
(\ref{Diophantine_condition}) and \eqref{condition_constant_b} yield
\[
	j' \ > \ v(j)\quad\textnormal{ with }\quad
	v(j) \ := \ a^{\frac{j-1}{dm}}+j.
\]
Hence, we obtain $\Leb_{\T^D}(\Omega_j)
\geq\Leb_{\T^D}(B_{r_{j+j_0-1}}(\tau_{j+j_0-1}))- \sum_{j'\geq
  v(j+j_0-1)}\Leb_{\T^D}(B_{r_{j'}}(\tau_{j'}))$, which is strictly positive if
$j_0\in\N$ is sufficiently large.  The corresponding subgraphs $\psi^j$ are
defined by restricting $\varphi^+$ to $\Omega_j$, i.e.\
$\psi^j:=\left.\varphi^+\right|_{\Omega_j}$.

\begin{prop} \label{dimensions_subgraphs} Let $T\in\mathcal{T}^*$.
  Then for all $j\in\set N_0$ the graph $\Psi^j$ is the image of a
  bi-Lipschitz continuous function $g_j:\Omega_j\to\Omega_j\times
  [0,1]$ and therefore $D_H(\Psi^j)=D$.  Further,
  $D_H(\Psi^\infty)\leq 1$.
\end{prop}
\proof
  Consider the maps $g_j:\Omega_j\to\Omega_j\times
  [0,1]:\theta\mapsto(\theta,\psi^j(\theta))$.  For all $j\in\set
  N_0\cup\{\infty\}$ we have $g_j(\Omega_j)=\Psi^j$ and
  $d_{\T^D\times[0,1]}(g_j(\theta),g_j(\theta'))\geq d(\theta,\theta')$ for all
  $\theta,\theta'\in\Omega_j$.  Further, for all $j\in\set N_0$ we have
\begin{align*}
		d_{\T^D\times[0,1]}(g_j(\theta),g_j(\theta'))
		\ \leq \ \left(1+K\alpha^{(j+j_0)m}\right)d(\theta,\theta'),
\end{align*}
for all $\theta,\theta'\in\Omega_j$.  This is true because Proposition
\ref{properties_approx_graphs} (iii) implies that
$\left.\varphi_n\right|_{\Omega_j}$ is Lipschitz continuous with Lipschitz
constant $K\alpha^{(j+j_0)m}$ independent of $n$, and since
$\psi^j=\lim_{n\to\infty}\left.\varphi_n\right|_{\Omega_j}$ we also get that
$\psi^j$ is Lipschitz continuous with the same constant.  This means that $g_j$
is bi-Lipschitz continuous for any $j\in\set N_0$, and therefore
$D_H(\Psi^j)=D_H(\Omega_j)$.  Hence, $D_H(\Psi^j)=D$ for all $j\in\set N_0$
because $0<\Leb(\Omega_j)<\infty$.

In order to complete the proof, we now show that $D_H(\Psi^\infty)\leq 1$.
Since $\Omega_\infty$ is a $\limsup$ set and for all $s>0$ we have
$\sum_{k=1}^\infty\diam(B_{r_k}(\tau_k))^s<\infty$, we get that
$D_H(\Omega_\infty)\leq s$ for all $s>0$, using Lemma
\ref{lemma_Hausdorff_dimension_limsup_set}.  Hence, $D_H(\Omega_\infty)=0$.
Furthermore, $\Psi^\infty\subset\Omega_\infty\times[0,1]$ and therefore
$D_H(\Psi^\infty)\leq D_H(\Omega_\infty)+D_B([0,1])=1$, applying Theorem
\ref{theorem_Hausdorff_dimension_product_sets}.
\qed\medskip

Since the Hausdorff dimension is countably stable, we immediately obtain

\begin{thm} \label{theorem_Hausdorff_dimension} Let
  $T\in\mathcal{T}^*$. Then the Hausdorff dimension of the upper
  bounding graph is $D$.
\end{thm}

It remains to determine the $D$-dimensional Hausdorff measure of $\Phi^+$.
\begin{prop} \label{proposition_Hausdorff_measure_finite} 
Let $T\in\mathcal{T}^*$ and $D>m^2\log(\alpha/a)$. Then
 the $D$-dimensional Hausdorff measure of $\Phi^+$ is finite.
\end{prop}
\proof Since $D_H(\Psi^\infty)\leq 1$, we have $\mathcal
  H^D(\Psi^\infty)=0$ for $D>1$.  Furthermore, we can consider the maps $g_j$
  from the last proposition as Lipschitz continuous maps from $\set R^D$ to
  $\set R^{D+1}$ and therefore we can use the Area formula (see for example
  \cite[Chapter 3]{EvansGariepy1992}) to deduce
	\begin{eqnarray*}
		\mathcal H^D(\Psi^j) \ &\leq& \ 
		\sqrt{1+(K\alpha^{(j+j_0)m+1})^2} 
		\ \Leb_{\set R^D}(B_{r_{j+j_0-1}}(\tau_{j+j_0-1}))\\
		\ &=& V_D\left(\frac{b}{2}\right)^D\sqrt{1+(K\alpha^{(j+j_0)m+1})^2} \ 
		a^{-\frac{D}{m}(j+j_0-2)}.
	\end{eqnarray*}
	When $D>m^2\log(\alpha/a)$ this implies that $\mathcal H^D(\Psi^j)$ is
        decaying exponential fast, and therefore $\mathcal H^D(\Phi^+) =
        \sum_{j=0}^\infty \mathcal H^D(\Psi^j)<\infty$.
\qed\medskip

\begin{prop} \label{proposition_Hausdorff_measure_infinite} Let
  $T\in\mathcal{T}^*$ and $D=1$. Then the one-dimensional Hausdorff
  measure of $\Phi^+$ is infinite.
\end{prop}
\proof
  We show that there exists an increasing sequence of integers
  $\ifolge{j_i}$ such that $\mathcal H^1(\Psi^{j_i})\geq c^+/6$.

  Suppose $j_1\ld j_N$ are given.  Our first goal is to find
  $j>j_N+j_0-1$ such that there exists a point $\tilde\theta^+\in
  B_{r_j}(\tau_j)$ with $\varphi_j(\tilde\theta^+)\geq 2c^+/3$.
  According to Remark \ref{r.0-lyap}, we can find a $\theta^+\in\T^1$
  with
  $\theta^+\notin\Omega_\infty':=\bigcap_{i=0}^\infty\bigcup_{k=i+1}^\infty
  B_{2r_k}(\tau_k)$ and $c^+:=\varphi^+(\theta^+)>0$. Since
  $\theta^+\notin\Omega_\infty'$, there exists $q\in\set N$ such that
  $\theta^+\notin\bigcup_{k=q}^\infty B_{2r_k}(\tau_k)$.  Now, we can
  choose $n>\max\{j_N+j_0-1, mq\}$ such that for all $j\geq n$
\begin{gather}
		\frac{1}{6}c^+\ \geq \ \frac{1}{1-\alpha^{-\lambda}}\alpha^{-\lambda j},
			\label{condition_distance_approx_graphs}\\
		v(j) \ \geq \ m(j+1)+1,\label{condition_stabilisation}\\
		a^{\frac{v(j)-1}{m}} \ \geq \ \frac{6 b}{c^+(1-a^{-1/m})}
			\left(1+K\alpha^{(j+1)m+1}\right)\label{condition_v_j}.
\end{gather}
Note that $B_{r_n}(\theta^+)\cap\bigcup_{k=q}^{n} B_{r_k}(\tau_k) =
\emptyset$, which means that there exists a neighbourhood of
$\theta^+$ where we can apply Proposition
\ref{properties_approx_graphs} (ii) to all points of this
neighbourhood.  Since $\varphi_n$ is continuous and
$\varphi_n(\theta^+)\geq \varphi^+(\theta^+)=c^+$, we can find
$\delta\leq r_n$ such that $\varphi_n(\theta)>5c^+/6$ for all
$\theta\in B_{\delta}(\theta^+)$.  Now, let $j\geq n$ be the first
time such that $B_\delta(\theta^+)\cap B_{r_j}(\tau_j)\neq\emptyset$.
Set $R:=B_\delta(\theta^+)\backslash B_{r_j}(\tau_j)\ne\emptyset$. Then for all
$\theta\in R$ we have $\theta\notin\bigcup_{k=q}^{n'} B_{r_k}(\tau_k)$
for all $n\leq n'\leq j$ and therefore
\[
        \sum\limits_{k=n}^{j-1}\alpha^{-\lambda k} \ \geq \ \varphi_n(\theta)-
        \varphi_j(\theta) \ > \ \frac{5c^+}{6}-\varphi_j(\theta),
\]
using $n\geq qm+1$ and Proposition \ref{properties_approx_graphs} (ii).  This
implies $\varphi_j(\theta)>2c^+/3$ for all $\theta\in R$, using
\eqref{condition_distance_approx_graphs}.  Since $\varphi_j$ is continuous, there
exists a $\tilde\theta^+\in B_{r_j}(\tau_j)$ such that
$\varphi_j(\tilde\theta^+)\geq 2c^+/3$.  Now, using Proposition
\ref{properties_approx_graphs} (i), we have that $\varphi_j$ is Lipschitz
continuous with Lipschitz constant $\beta\alpha^j$ and therefore there exists an
interval $I\ssq B_{r_j}(\tau_j)$ such that $\varphi_j$ is greater than $c^+/2$
on $I$ and
	\[
		\Leb_{\T^1}(I) \ \geq \ \frac{c^+}{6\beta\alpha^j}.
	\]
	Because of \eqref{condition_v_j}, we have that
	$\Leb_{\T^1}(I\backslash\bigcup_{k=j+1}^{\infty} B_{r_k}(\tau_k))>0$ (note
	that $\beta< K$).
	Hence, using \eqref{condition_stabilisation} plus Proposition
	\ref{properties_approx_graphs} (ii) and
	\eqref{condition_distance_approx_graphs} again, there exists $\theta\in
	I\backslash\bigcup_{k=j+1}^{\infty} B_{r_k}(\tau_k)\subset\Omega_{j_{N+1}}$
	such that $\psi^{j_{N+1}}(\theta)\geq c^+/3$, where $j_{N+1}:=j-j_0+1$.
	Finally, the application of \eqref{condition_v_j} yields
	\begin{eqnarray*}
          \mathcal H^1(\Psi^{j_{N+1}}) \  &\geq& \
          \mathcal H^1(\psi^{j_{N+1}}(\Omega_{j_{N+1}}))\\  
          \  &\geq& \  \frac{c^+}{3}-\left(1+K\alpha^{(j+1)m+1}\right)
          \Leb_{\T^1}\left(\bigcup\limits_{k=j+1}^{\infty} B_{r_k}(\tau_k)\right)
          \ \geq \ \frac{c^+}{6} \ . 
	\end{eqnarray*}
\qed\medskip

We turn to the question of rectifiability.  Note that by definition
$\mu_{\varphi^+}$ is absolutely continuous with respect to
$\left.\mathcal H^D\right|_{\Phi^+}$.

\begin{thm} \label{theorem_rectifiability} Let $T\in\mathcal{T}^*$.
  Then $\mu_{\varphi^+}$ is $D$-rectifiable and
  $d_{\mu_{\varphi^+}}=D_1(\mu_{\varphi^+})=D$.
\end{thm}
\proof
	Observe that $\mu_{\varphi^+}(\Psi^\infty)=0$.
	Therefore, $\mu_{\varphi^+}$ is also absolutely continuous with respect
	to $\left.\mathcal H^D\right|_{\Phi^+\backslash\Psi^\infty}$ and
	$\Phi^+\backslash\Psi^\infty=\bigcup_{j=0}^\infty\Psi^j$
	is countably $D$-rectifiable, according to Proposition \ref{dimensions_subgraphs}.
	That means $\mu_{\varphi^+}$ is $D$-rectifiable.
	Now, use Corollary \ref{dimensions_rectifiable_measure} to obtain the
	dimensional results for $\mu_{\varphi^+}$.
\qed\medskip

Note that for $D\geq 2$ we have $\mathcal H^D(\Psi^\infty)=0$, such that
$\Phi^+$ is countably $D$-rectifiable.  The question whether $\Phi^+$ is
countably $1$-rectifiable for $D=1$ remains open.

We can now apply the above results to the family $F_\kappa$ defined in
Example~\ref{example_higher_dimensional_case} to obtain the following
corollary, which contains Theorem~\ref{t.hausdorff-dim} and 
\ref{t.pointwise} and Corollary~\ref{t.info-dim} as a special case.
\begin{cor}
  \label{c.final} Let $F_\kappa$ be defined by (\ref{e.F_kappa}). Then
  there exists a $\kappa_0=\kappa_0(c,d,D)$ such that for all $\kappa\geq
  \kappa_0$ 
  \begin{itemize}
  \item the upper bounding graph $\Phi^+$ of $F_\kappa$ has
  Hausdorff dimension $D$;
\item the $D$-dimensional Hausdorff measure of $\Phi^+$ is infinite if
  $D=1$ and finite for $D$ sufficiently large;
  \item  $\mu_{\varphi^+}$ is exact dimensional with
  pointwise dimension $D$;
\item the information dimension of $\mu_{\varphi^+}$ is $D$;
\item $\mu_{\varphi^+}$ is $D$-rectifiable.
\end{itemize}
\end{cor}

Finally, we close by addressing a further obvious question in our context,
namely that of the size of the set of {\em `pinched points'} where the upper
bounding graph $\varphi^+$ equals zero. Given $T\in\mathcal{T}$, let
\[
\mathcal{P} \ := \ \left\{ \theta \in \T^D \mid \varphi^+(\theta)=0 \right\} \ .
\]
Then $\mathcal{P}$ is residual in the sense of Baire
\cite{Keller1996}, and therefore its box dimension and its packing
dimension are $D$. However, from the point of view of Hausdorff
dimension, $\mathcal{P}$ turns out to be small. 

\begin{prop} \label{p.pinchedset}
  Let $T\in\mathcal{T}^*$. Then
  \[
  \mathcal{P} \ \ssq \ \Omega_\infty \cup \left\{\omega^n(\theta_*)\mid
  n\in\N\right\} \ ,
  \]
  where $\Omega_\infty$ is the set defined in (\ref{e.omega_infty}).
  In particular, $D_H(\mathcal{P}) = 0$.
\end{prop}
\proof Suppose $\theta\notin \Omega_\infty \cup \left\{\omega^n(\theta_*)\mid
  n\in\N\right\}$. Let $q\in\N$ be such that $\theta\notin \bigcup_{j=q}^\infty
B_{r_j}(\tau_j)$ and fix any $t\geq mq$. Let 
\[
\eps \ := \ \min_{k=1}^t T^k_{\omega^{-k}(\theta)}(L_0) \ .
\]
Note that since $\theta\notin \{\omega^n(\theta_*)\mid n\in\N\}$ we have $\eps >
0$. Now, for any $n>t$ Lemma~\ref{main_property_of_s_n} implies that
$s^n_{n-t}(\theta) \leq 11t/m\leq t/2$. In particular, there exists
$l\in\{n-t\ld n-1\}$ such that $x_l=T_{\omega^{-n}(\theta)}^l(1) \geq L_0$. Hence,
\[
\varphi_n(\theta) \ = \ T^{n-l}_{\omega^{-(n-l)}(\theta)}(x_l) \ \geq \ \eps \ .
\]
Since this holds for all $n>t$, we obtain $\varphi^+(\theta)\geq \eps$ and
thus $\theta\notin \mathcal{P}$ as required. The statement on the Hausdorff
dimension then follows from Lemma
\ref{lemma_Hausdorff_dimension_limsup_set}. \qed

\end{document}